\documentclass[11pt]{article}
\usepackage{graphicx,float,xcolor} 
\usepackage{amsthm,amsfonts,amsmath}
\usepackage{amssymb,bm,bbm,authblk}
\usepackage{cite}
\usepackage{hyperref,sectsty,titlesec}

\usepackage[top=1in, bottom=1in, left=1in, right=1in]{geometry}
\sectionfont{\fontsize{12}{14}\selectfont}
\titlespacing{\section}{0pt}{\parskip}{\parskip}

\title{What Do Bouncing Balls Tell Us About the Universe? A Journey into Billiard Systems}
\author[1]{Weiqi Chu}
\author[1]{Matthew Dobson}
\affil[1]{Department of Mathematics and Statistics, University of Massachusetts Amherst}
\date{}

\begin{document}

\maketitle

\section*{Abstract}
Have you ever played or watched a game of pool? If so, you have already seen a billiard system in action. In mathematics and physics, a billiard system describes a ball that moves in straight lines and bounces off walls. Despite these simple rules, billiard systems can produce remarkably rich behaviors: some table shapes generate regular, periodic patterns, while others give rise to complete chaos. Scientists also study what happens when we shrink the ball down to the size of an electron to a world where quantum effects take over and the familiar reflection rules no longer apply. In this article, we discuss billiard systems in their many forms and show how such a simple setup can reveal fundamental insights into the behavior of nature at both classical and quantum scales.

\section*{What are billiard systems?}
A billiard system is a ball moving in straight lines and bouncing off walls~\cite{sinai2004billiard}. The walls of the table can take many shapes: a rectangle, a circle, or a complex table with a mix of straight and curved walls. When the ball hits a wall, it bounces off at the same angle it arrived, going the other way, just like light bouncing off a mirror. The incoming angle is known as the angle of incidence, and the outgoing angle is known as the angle of reflection. In what is known as the Law of Reflection, the angle of incidence equals the angle of reflection.

\begin{figure}[htp]
    \centering
    \includegraphics[width=0.9\linewidth]{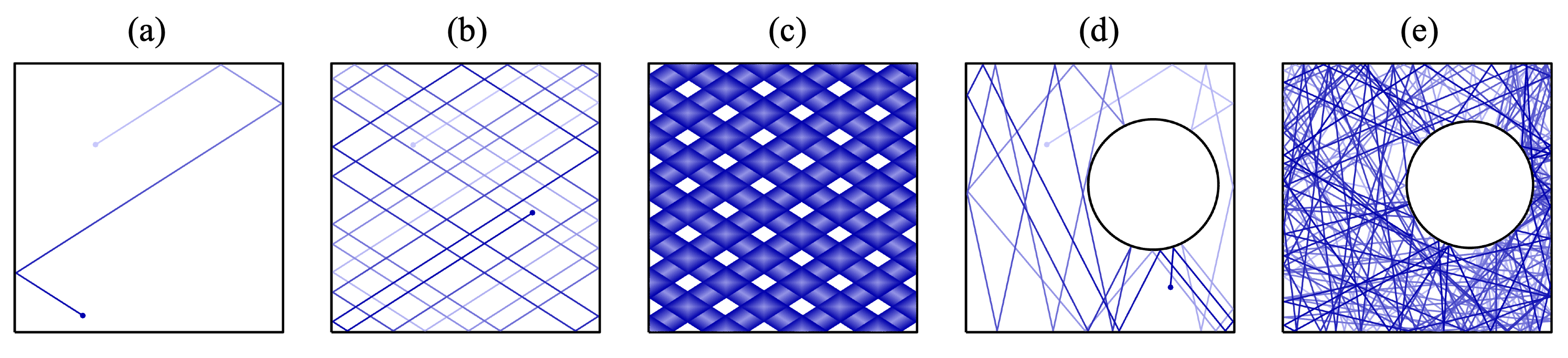} 
    \caption{(a--c) A square billiard. The trajectory of a ball on a square table after 3, 30, and 300 bounces. By the time the ball bounces 300 times, the paths overlap and create an easily recognizable pattern. Moreover, every collision with the top wall occurs at the same angle, no matter how many bounces take place. (d, e) A Sinai billiard. We add a circular obstacle inside the square perimeter, so the ball reflects off both the walls and the circle. The figures show the trajectory of a ball after 30 and 300 bounces, and the trajectory no longer forms a regular pattern but instead appears disordered. The reflection angles against the top wall are no longer the same for each collision.}
    \label{fig:examples}
\end{figure}

Even with such a simple setup, the ball's path can reveal surprising patterns and behaviors. Depending on the shape of the walls, the ball can either trace out recognizable patterns or move chaotically. Figure~\ref{fig:examples} shows two billiard systems: a square billiard and a Sinai billiard. A Sinai billiard is a square-shaped table with a circular obstacle placed inside. In both cases, the ball starts at the same position and with the same shooting angle, yet the paths are very different.

So why do scientists and mathematicians care about these bouncing balls? Billiard systems can reveal the rules that govern motion beyond just a pool table. By studying the trajectory of a ball reflecting off boundaries, researchers can identify what separates order from chaos and connect these ideas to deep areas of research. Beyond theory, billiards have practical applications in the real world, like light rays in optical cavities and sound waves in concert halls, and can even show us the optimal trajectories of robots moving around. In this article, we will trace the path of a bouncing ball and discover how such a simple concept unlocks profound insights in physics and mathematics.

\section*{Integrable and chaotic billiards}
The properties of billiard systems depend on the shape of the walls and where we start the ball.  The ball's behavior can range from completely regular to complete chaos~\cite{chernov2006chaotic}. Some shapes produce very regular motion where a ball bounces in repeating patterns or even retraces the same loop with reflection angles preserved throughout the trajectory. Billiards with these properties are called integrable billiards. 
In contrast, if you observe disordered paths, the systems are likely chaotic billiards. Between completely integrable and fully chaotic billiards, there are mixed billiards, where some initial conditions lead to regular motion while others produce chaotic behavior.

\begin{figure}[htp]
    \centering
    \includegraphics[width=0.75\linewidth]{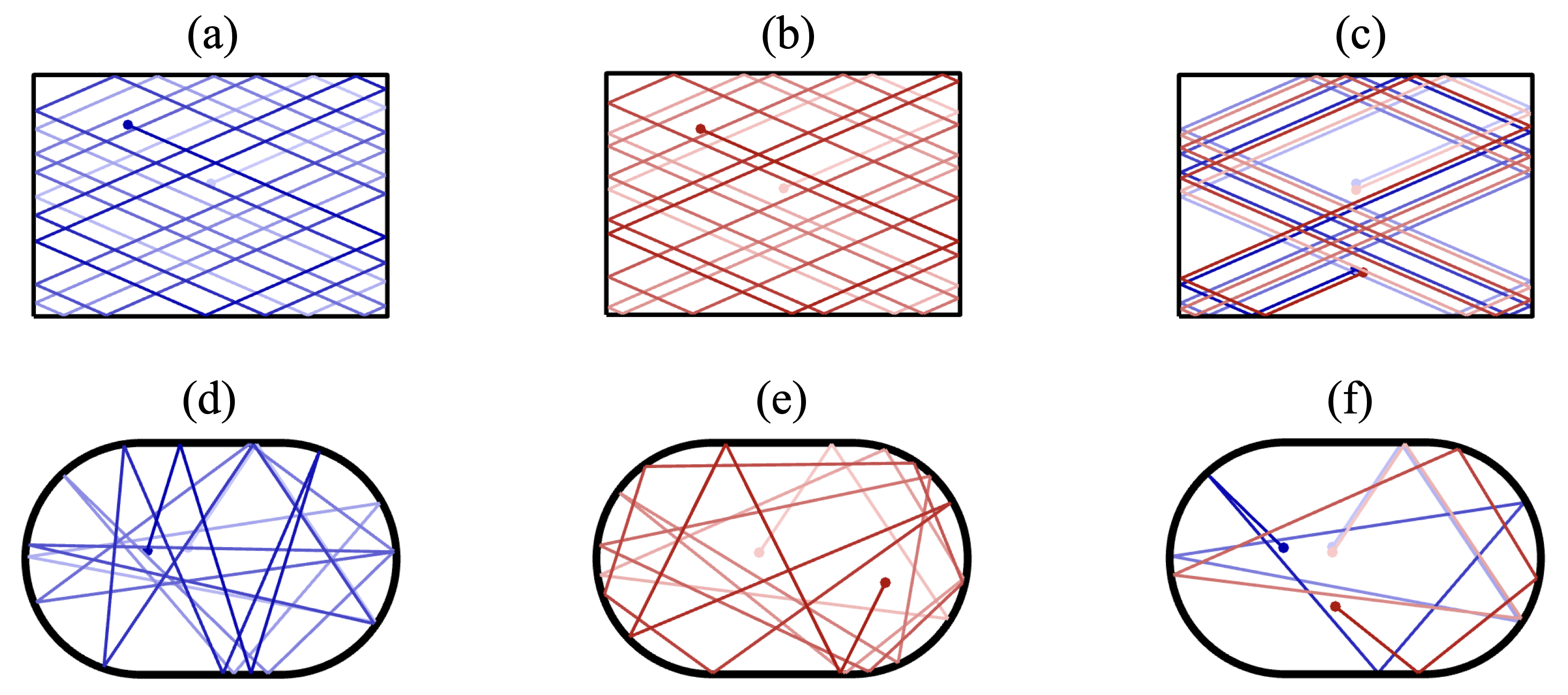}
    \caption{(a--c) In a square billiard, two balls start at almost the same position and their trajectories are shown in (a) and (b). We plot the paths in the same figure (c), and the two balls remain close to each other after 30 bounces. (d--f) In a Bunimovich stadium billiard, two balls spread apart and are no longer close after 6 bounces, as shown in (f). The entire trajectories of these two balls are also very different, as shown in (d) and (e).}
    \label{fig: chaotic}
\end{figure}

Let us look at two tables and see integrability and chaos.
Figure~\ref{fig: chaotic}(a--c) show an example of the square billiard. Two balls start from almost the same position in the same direction, and their paths remain close together throughout their trajectories. Even after many bounces, their paths stay nearby, and the reflection angles remain the same. The square billiard is integrable.  
Figure~\ref{fig: chaotic}(d,e) show the Bunimovich stadium billiard, which has a flat top and bottom like the square table, but has semicircles on each end. Two balls again start from almost the same position, with the red trajectory initially overlapping the blue one. Once they strike the circular boundary, their angles begin to split, and their paths separate more and more with each bounce. In this case, we say the Bunimovich stadium billiard is \emph{sensitive to its initial conditions}, which is a key feature of chaotic billiards. This initial condition sensitivity is also related to the well-known butterfly effect, the idea that a butterfly flapping its wings in one part of the world could, through a chain of tiny influences, lead to a tornado elsewhere. In chaotic billiards, tiny changes in initial conditions grow quickly and can compound into large differences in the long run, as the two balls in the Bunimovich stadium billiard.

Integrable billiards are not sensitive to initial conditions and also lack a property called \emph{mixing}. The mixing property describes how motion is scrambled over time, like stirring cream into coffee until it is perfectly blended. It requires a single trajectory to wander through all regions rather than remain confined to a repeating pattern. In an integrable billiard, the ball follows regular, repeating paths as shown in Figure~\ref{fig:examples}(c). In contrast, in a chaotic billiard, a ball's path spreads into many parts of the table, eventually covering the entire space in a seemingly random way. Chaotic billiards can have this property, while integrable billiards do not. In a mixing system, the connection to previous trajectories gradually fades. As the system evolves, its current state becomes less and less dependent on how it started. Eventually, every region is visited, and the overall pattern appears random and uniform, regardless of the initial conditions. This ``forgetting the past'' is what allows a mixing system to fully explore its space.

\section*{Quantum billiards}
So far, all the billiard systems we have discussed belong to the classical domain, where the balls are relatively large and visible to the naked eye. But what happens when we shrink down to the tiny world of nucleons and electrons and let these particles bounce off walls? At such small scales, usually around a billionth of a meter (nanometers) or smaller, classical physics, such as Newton's laws, no longer applies. Instead, we must use the rules of quantum mechanics, a fundamental branch of physics that describes how matter and energy behave at the subatomic level. Quantum billiards explore how particles like electrons behave when trapped inside microscopic cavities, and their studies blend the rules of quantum mechanics with the dynamics of classical billiards.

In the quantum world, particles behave very differently from billiard balls on a pool table. In classical physics, we can watch balls roll, predict their paths, and know exactly where they are and how fast they are moving. But in the strange world of quantum particles, things work differently. We cannot know the location at the same time as the speed and direction of a quantum particle, such as an electron. 
Because of this, a particle is not like a tiny marble rolling on a table. Instead, it is more like a cloud of possibility. Even though we cannot pin down the precise position of an electron, we can describe the places where it is most likely to appear. To do this, scientists use something called a wave function. A wave function acts like a special kind of wave --- not a wave of water, but a wave of possibility. It does not tell us exactly where the particle is. Instead, it shows the chances of finding the particle in certain areas. Where the wave is strong, the electron is more likely to be found. Where the wave is weak, it is less likely. 
Scientists use wave functions to describe how electrons move. These wave functions can travel, bounce, and reflect off boundaries much like a ball on a table. This idea forms the basis of quantum billiards.

\begin{figure}[htp]
    \centering
    \includegraphics[width=0.9\linewidth]{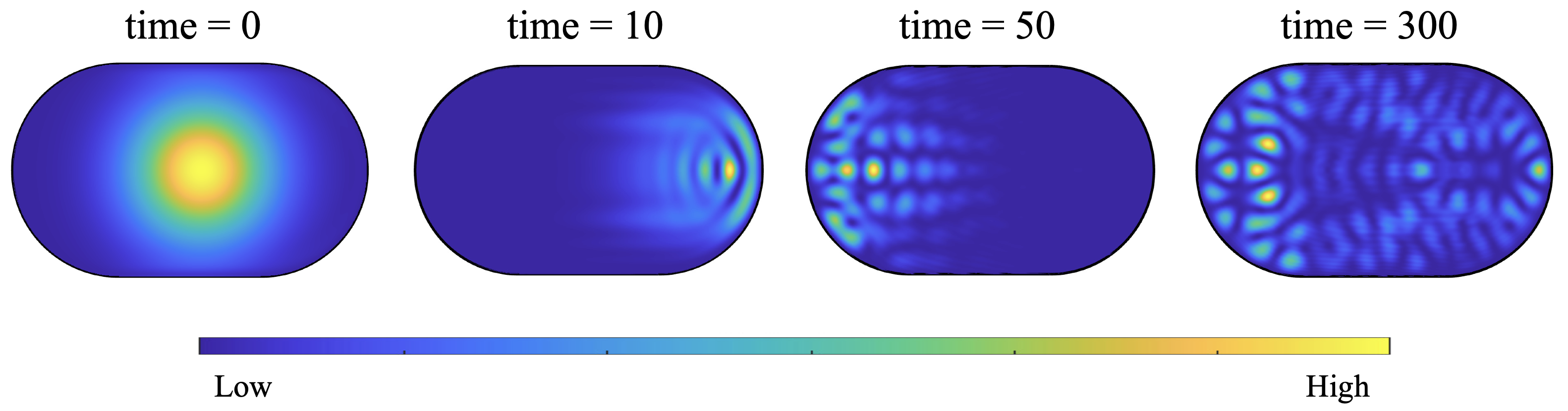}
    \caption{The figures show the likely location of an electron at any given moment, through color. At the beginning, the electron is likely in the center of the pool and is launched to the right. The wave function bounces off the curved wall on the right and reflects back, forming a ripple. When time is around $50$, the wave function reaches the left wall and bounces back again. After a few bounces, the wave function spreads out, and the electron can potentially be found anywhere within the stadium.}
    \label{fig: quantum}
\end{figure}

In a quantum billiard, wave functions move and bounce off the walls of the table, much like ripples on a pond after tossing in a pebble. Figure~\ref{fig: quantum} shows an example in a stadium-shaped billiard and shows the wave functions in color. The bright and dark spots (yellow and blue regions) mark the areas where the electron is more or less likely to be found. The wave function concentrates in the center originally, and when it is pushed to the right, it travels across, bounces off the walls, and reflects just like ripples on water. 

In the last section, we saw that in chaotic billiards (such as Figure~\ref{fig: chaotic}(f)), even the tiniest change at the start can lead to very different paths for the balls, while in integrable billiards (such as Figure~\ref{fig: chaotic}(c)), the paths stay close to each other. In quantum billiards, things are trickier. We cannot track the exact path of tiny particles like electrons, so we need different tools to figure out whether the system is chaotic or not. This is not an easy task. Indeed, classical and quantum billiards can show very different chaotic behavior, even if they are played on the same table shape. The study of this connection is called quantum chaos~\cite{jensen1992quantum}. It looks at how quantum systems behave when their classical versions are chaotic, and whether the rules of quantum theory can describe classically chaotic systems.

\section*{Final bounce}

From bouncing balls to wave functions, billiard systems provide a compelling window into both classical and quantum physics. We began with the simple scenario of a ball bouncing around a table and examined how different table shapes can result in fundamentally different behaviors, distinguishing integrable and chaotic systems. We explored how sensitive billiard systems are to initial conditions and how mixing can cause chaotic systems to eventually end up the same. 
Finally, we zoomed in to the scale of electrons and observed how wave functions act and bounce within a space, bringing the study of billiard dynamics to the quantum world.

Billiard systems extend far beyond the basic scenarios that we have discussed so far. For example, what happens if a small, random disturbance is applied each time the ball strikes a wall, or if the ball loses energy after each bounce and gradually slows down, or even if the ball is allowed to spin as it moves? These variations enrich the dynamics and bring the models closer to the real world, like fluids moving through pipes and pools, light reflecting inside crystals, or even a spinning dodge ball bouncing around. There are many questions about billiards systems that mathematicians are actively working on~\cite{bialy2022open}.

Billiards begin with a simple rule of motion, yet they open a gateway to exploring complex and fundamental questions in science. The simple act of bouncing a ball provides insight into the concepts of order, chaos, wave phenomena, and the strange behavior of the quantum world. Perhaps the next game of pool will spark a new scientific discovery.

\bibliographystyle{abbrv}
\bibliography{references}

\section*{Acknowledgment}

The authors thank Mathew A. Gan for his helpful feedback. 
W. Chu acknowledges support from NSF under grant DMS-2309814.
The authors used ChatGPT to help refine the writing and simplify the content, making it more approachable for a young audience.

\section*{Glossary}

\textbf{Initial conditions}: Initial conditions describe how something begins. For example, if you roll a ball on a billiard table, its starting position, direction, and speed are the initial conditions.

\textbf{Chaos}: Chaos means a state that is very unordered and mixed up. In a chaotic system, things get tangled and spread out so much that small differences at the start can lead to completely different outcomes.

\textbf{Nucleon}: A nucleon is a subatomic particle in the nucleus of an atom, either a proton or a neutron. 

\textbf{Wave function}: A wave function is a mathematical description of the state of an electron. Where the magnitude of the wave function is large, the electron is more likely to be found there.

\section*{Author introduction}
Weiqi Chu is an Assistant Professor in the Department of Mathematics and Statistics at the University of Massachusetts Amherst. Her research focuses on scientific computing and network science. Outside of mathematics, Weiqi enjoys reading science fiction, playing strategy card games, and exploring local cuisines while traveling. 

Matthew Dobson is an Associate Professor in the Department of Mathematics and Statistics at the University of Massachusetts Amherst. His research is on how to simulate the behavior of matter at the molecular scale. 

\end{document}